\documentclass[12pt]{article}
\usepackage{amssymb,latexsym,amsmath}

\newtheorem{theorem}{Theorem}[section]

\newtheorem{corollary}[theorem]{Corollary}
\newtheorem{definition}[theorem]{Definition}
\newtheorem{conjecture}{Conjecture}

\newenvironment{proof}{\paragraph{\it Proof.}}{$\square$\vskip0.4cm}
\newenvironment{remark}{\paragraph{\it Remark.}}{\vskip0.4cm}

\def\endproof{\hfill$\square$\medskip}

\def\be{\begin{eqnarray}}
\def\ee{\end{eqnarray}}
\def\bes{\begin{eqnarray*}}
\def\ees{\end{eqnarray*}}
\def\A{\mathcal{A}}

\def\B{\mathcal{B}}

\def\C{\mathbb{C}}
\def\M{\mathcal{M}}

\def\R{\mathbb{R}}
\def\Q{\mathbb{Q}}

\def\Z{\mathbb{Z}}
\def\N{\mathbb{N}}
\def\K{\mathbb{K}}

\begin{document}
\thispagestyle{empty}

\title{Quaternionic Geometry of Matroids}

\author{
 Tam\'as Hausel
\\ {\it Department of Mathematics}
\\ {\it University of Texas at Austin}
\\ {\it Austin TX 78712, USA}
\\{\tt hausel@math.utexas.edu}}

\maketitle
\begin{abstract} Building on a recent paper \cite{hausel-sturmfels}, 
here we argue that the combinatorics of 
matroids is intimately related to the geometry and topology of toric 
hyperk\"ahler varieties. We show that just like toric varieties occupy
a central role in  Stanley's proof for the necessity of McMullen's
conjecture (or $g$-inequalities) about the classification of face
vectors of simplicial polytopes, the topology of toric hyperk\"ahler
varieties leads to new restrictions on face vectors of matroid
complexes. Namely in this paper we will give two proofs that the
injectivity part of the Hard Lefschetz theorem survives for toric
hyperk\"ahler varieties. We explain how this implies the
$g$-inequalities for rationally representable matroids. We show how
the geometrical intuition in the first proof, coupled with results of
Chari \cite{chari}, leads to a proof of the $g$-inequalities for
general matroid complexes, which is a recent result of Swartz
\cite{swartz}. 
The geometrical idea in the second proof 
will show that a 
pure $O$-sequence should satisfy the $g$-inequalities, 
thus showing that our result is in fact a consequence of a long-standing conjecture of Stanley.

\end{abstract}

\section{Introduction}

McMullen \cite{mcmullen} conjectured in 1971 that the face vector\footnote{$f_i$ is the number of
$i$-dimensional faces}  
$(f_0,\dots,f_{k-1})$ of a $k$-dimensional simplicial
polytope $P\subset \R^k$ should satisfy, what we call, {\em $g$-inequalities}:
\begin{eqnarray}
\begin{array}{l} g_i\geq 0, \mbox{for } 1\leq i\leq \lfloor {\frac{k}{2}} \rfloor,\\ 
\mbox{and, if one writes }\\
g_i=\binom{n_i}{i}+\binom{n_{i-1}}{i-1}+\dots + \binom{n_r}{r},\\
\mbox{with } n_i>n_{i-1}>\dots>n_r\geq r\geq 1, \mbox{ then } \\
 g_{i+1}\leq \binom{n_{i}+1}{i+1}+\binom{n_{i-1}+1}{i}+\dots + \binom{n_r+1}{r+1}\\
\mbox{ for } 1\leq i<\lfloor {\frac{k}{2}} \rfloor,
\label{ginequalities}
\end{array}
\end{eqnarray}

where $$g_i=h_i-h_{i-1}$$ and \begin{eqnarray}\label{hi}h_i=\sum_{j=i}^k(-1)^{j-i}\binom{j}{i}f_{k-j-1}.\end{eqnarray}

Stanley \cite{stanley-toric} in 1980 proved this conjecture using toric 
varieties. In a nutshell his argument goes as follows. First one perturbs the vertices of 
$P$ a little bit so that $P$ becomes a rational polytope. Because $P$
is simplicial 
this does not change the face vector of $P$. The next step is to take
the corresponding $k$-dimensional toric orbifold $X(\Delta_P)$, where
$\Delta_P$ is the fan of cones over the faces of $P$. It is a well-known fact (see e.g. \cite{fulton}) that the  $i$'th $h$-number $h_i=b_{2i}(X(\Delta_P))$ agrees with the $2i$th Betti number of $X(\Delta_P)$. 
Now $X(\Delta_P)$  has an ample class $\omega\in H^2(X(\Delta_P), \C)$, which induces
a map $$L:H^*(X(\Delta_P),\C)\to H^*(X(\Delta_P),\C),$$ by multiplication with $\omega$. Using the injectivity part of the Hard Lefschetz theorem (see e.g. \cite{griffiths-harris}), 
which implies that $L$ is an injection below 
degree $k$,  we get that the degree $2i$th part of the  graded algebra $H^*(X(\Delta_P),\C)/(\ker (L))$ has dimension 
\begin{eqnarray}\dim(H^{2i}(X(\Delta_P),\C)/(\ker (L)))=h_i-h_{i-1}=g_i\label{dimension}
\end{eqnarray} for $2i<k$. Since $H^*(X(\Delta_P),\C)$ is generated by $H^2(X(\Delta_P),\C)$ we also get that the algebra   $H^*(X(\Delta_P),\C)/(\ker (L))$ is generated in degree $2$. 
Now, using (\ref{dimension}), a well-known theorem of Macaulay (see
e.g. \cite[Theorem II.2.3]{stanley-book}) proves the g-inequalities (\ref{ginequalities}). See \cite{fulton} or \cite{stanley-book} for more details.  

Our starting point is the observation \cite[Corollary 1.2]{hausel-sturmfels} that the $h$-vectors of a rationally representable matroid 
$\M_\B$ agree $h_i(\M_\B)=b_{2i}(Y(A,\theta))$ with the  Betti numbers
of a toric hyperk\"ahler variety $Y(A,\theta)$, 
for a generic choice of 
$\theta$, where the toric hyperk\"ahler variety can be considered as a quaternionic analogue of a toric variety. 
Therefore any restriction on the cohomology of a toric hyperk\"ahler variety will yield restrictions on the face vectors of rationally representable matroid 
complexes and vice versa any known restriction on the face vectors of
(rationally representable) matroids yields cohomological restrictions on toric
hyperk\"ahler varieties. This two-way relationship between these two seemingly unrelated subjects, hyperk\"ahler geometry on one hand and combinatorics of matroids on the other, is 
what we call ``Quaternionic geometry of matroids''. 
A relationship of this flavor is exploited in a recent paper
by Swartz and the author \cite{hausel-swartz}. There the combinatorics of affine hyperplane arrangements yields the existence of many $L^2$ harmonic forms on the 
corresponding toric hyperk\"ahler manifold, in harmony with conjectures by physicists in string theory. For details see the paper \cite{hausel-swartz}. 

In the present paper our purpose is to use intuition arising from the study of the geometry of toric hyperk\"ahler varieties to prove results in the combinatorics of matroids. 
Namely we will proceed as follows: In the next Section~\ref{simplicial} and
Section~\ref{toric} we recall some basic notations and results
from  \cite{stanley-book}  and from \cite{hausel-sturmfels}. Then we go on and in Section~\ref{injective} give two 
different proofs for the injectivity part of the Hard Lefschetz theorem for toric hyperk\"ahler varieties. The second one is basically taken from \cite[Theorem 7.4]{stanley-book}, 
while the first proof could be easily generalized for other 
similar hyperk\"ahler manifolds, like for example Nakajima's quiver
varieties \cite{nakajima} or Hitchin's moduli of Higgs bundles \cite{hitchin}. In Section~\ref{g} then we explain how the geometric idea in the first proof can be generalized for 
any matroid complexes, a result recently proven by Swartz in \cite{swartz}. What we show is that the geometrical structure for the first proof is provided for general matroids by Chari's decomposition theorem \cite{chari}. 
In fact this proof is 
similar to Swartz's original proof in \cite{swartz}. We conclude our paper by showing that the geometric structure which yielded the second proof of the injective Hard Lefschetz theorem is 
present for pure 
$O$-sequences. This way we find that the $g$-inequalities we proved in the previous section are in fact a consequence of a long standing conjecture of Stanley \cite{stanley-conjecture}. This last result
is a strengthening of a result of Hibi in \cite{hibi}.
 
\paragraph{\bf Acknowledgement.} This paper grew out from a project started with Bernd Sturmfels in \cite{hausel-sturmfels}. Conversations with Edward Swartz were also useful. Financial
support was provided by a Miller Research Fellowship at the University
of California at Berkeley, and by NSF grants DMS-0072675 and DMS-0305505. 

\section{Simplicial and matroid complexes}
\label{simplicial}

We collect here some basic definitions and results on simplicial
complexes and in particular matroid complexes from
\cite{stanley-book}.

A simplicial complex $\Sigma$ on a finite set $V=\{1,\dots,n\}$ is a set of subsets
of $V$, i.e. $\Sigma\subset 2^V$, such that $\{x\}\in \Sigma$ for any
$x\in V$ and $F\in \Sigma$ and $F^\prime\subset F$ implies $F^\prime
\in \Sigma$. We call $F\in \Sigma$ a {\em face} of $\Sigma$, the
{\em dimension} of the face is one less than its size. The dimension
of $\Sigma$ is then the maximum dimension of its faces, while its {\em
  rank} is $1$ more. A {\em facet} is a
face of maximal dimension. A simplicial complex is called {\em pure}
if its maximal faces are all facets. 
The {\em $f$-vector} of a rank-$k$ 
simplicial complex is $(f_0,f_1,\dots,f_{k-1})$, where $f_i$ is the number
of $i$-dimensional faces in $\Sigma$. The {\em $h$-vector} of the
simplicial complex is $(h_0,\dots,h_{k})$ given by (\ref{hi}).

Define the {\em Stanley-Raisner ring} of a rank-$k$ 
simplicial complex $\Sigma$ 
as a graded ring given by:
$$\C[\Sigma]=\C[x_1,\dots,x_n]/\langle x_F=\prod_{i\in F} x_i |
F\notin \Sigma\rangle.$$ All our simplicial complexes in this paper
will be {\em Cohen-Macaulay}, which will imply that we will always
have a {\em linear system of operators} or {\em l.s.o.p} for short,
which is a sequence $(\theta)=(\theta_1,\dots,\theta_{k})$ of linear
combinations of the $x_i$, such that the graded ring
$$\C[\Sigma]/(\theta):=\C[\Sigma]/(\theta_1\C[\Sigma]+\dots+\theta_k\C[\Sigma])$$ is finite
dimensional as a vector space over $\C$ and that the $h$-numbers 
$h_i(\Sigma)=(\C[\Sigma]/(\theta))_i$ agree with the dimension of the 
corresponding graded
piece of $\C[\Sigma]/(\theta)$. 

We will use the following operation on simplicial complexes in
Section~\ref{g}. 
Given two simplical complexes $\Sigma$
with vertex set $V$ and $\Theta$ with vertex set $U$ we define their
{\em poset-theoretic} product $\Sigma\times\Theta$ as a simplicial
complex with vertex set $U\cup V$ and all faces of the form $F\cup
F^\prime$ where $F\in \Sigma$ and $F^\prime\in \Theta$. The
poset-theoretic product has the advantage that it behaves nicely after
taking the corresponding Stanley-Raisner rings: $\C[\Sigma\times
  \Theta]\cong \C[\Sigma]\otimes \C[\Theta].$

For examples of (Cohen-Macaulay) 
simplicial complexes we mention the boundary complex of 
a simplicial convex polytope, which was mentioned in the
introduction. Another class for interest for us are {\em matroid
  complexes} or simply just {\em matroids}. 
A matroid complex $\M$ is a simplicial complex on a vertex
set $V$ such that for every $W\subset V$ the induced subcomplex
$\M_W=\{F\in \M: F\subset W\}$ is pure. The {\em rank} of the
matroid is $1$ more than its dimension. A vertex $i\in V$ is a coloop of $\M$
if $\M_{V\setminus i}$ has rank smaller than the rank of $\M$. 

The motivating example of a
matroid complex $\M_\B$ on vertex set $V=\{1,\dots,n\}$ is obtained from a vector configuration
$\B=(b_1,\dots,b_n)\in \K^k$ in a $k$-dimensional vector space over a
field $\K$, defined by $F\in \M$ iff $\{b_i\}_{i\in F}$ is linearly
independent. Such a matroid is called {\em representable} over $\K$. For
example, if $\K=\Q$ then we call the matroid $\M$  rationally representable.   
 
For more details on these definitions consult \cite{stanley-book}, the
poset-theoretic product was used in \cite{chari}.

\section{Toric hyperk\"ahler varieties}
\label{toric}

Here we collect notation and terminology from \cite{hausel-sturmfels}
which we will need in the present paper. For more details see
\cite{hausel-sturmfels}. 

Let $A=[a_1,\dots,a_n]$ be a $d \times n$-integer matrix  whose
$d \times d$-minors are relatively prime. We choose
an $n \times (n \! - \! d)$-matrix  $B=[b_1,\dots,b_n]^T$ which makes
 the following sequence exact:
\begin{eqnarray*}0 \, \longrightarrow \, \Z^{n-d}
\, \stackrel{B}{\longrightarrow}  \, \Z^n \,
\stackrel{A}{\longrightarrow} \, \Z^d \,
 \longrightarrow \, 0.\label{ses}\end{eqnarray*} Taking $\theta\in
\N\A$, where $\A := \{a_1,\ldots,a_n\} $ is a vector configuration in
$\Z^d$, \cite{hausel-sturmfels} constructs a quasi-projective
variety $Y(A,\theta)$ (which sometimes we abbreviate as $Y$), 
called a {\em toric hyperk\"ahler variety}. (This construction is 
an algebraic geometric version of the original construction of 
Bielawski and Dancer in \cite{bielawski-dancer}.) By
\cite[Proposition 6.2]{stanley-book} if $\theta\in \N
\A $ is generic $Y(A,\theta)$ is an orbifold, while 
if, in addition, $A$ is unimodular then $Y(A,\theta)$ is a smooth variety. 

The topology of $Y(A,\theta)$ is governed by an affine hyperplane
arrangement denoted by ${\mathcal H}(B,\psi)$ of $n$ planes in
$\R^{n-d}$. For example a key result in \cite[Corollary 6.6]{stanley-book} claims
that the $h$-numbers of the matroid of the vector configuration
$\B=\{b_1,\dots,b_n\}$ 
agree with the Betti
numbers of $Y$: $$h_i(\M_\B)=b_{2i}(Y(A,\theta)).$$ In the next section
we will make use of a projective subvariety $C(A,\theta)$ of
$Y(A,\theta)$, which is called the {\em core} of $Y(A,\theta)$. It is a
reducible variety whose components are projective toric varieties,
corresponding to top dimensional bounded regions in ${\mathcal
  H}(B,\psi)$. If the matroid of $\B$ is coloop-free than the core is a
middle and pure dimensional projective subvariety of $Y(A,\theta)$. 

Finally we need to mention a result from \cite{harada-proudfoot}. They
construct and study a certain residual $U(1)$-action on $Y(A,\theta)$,
which comes from an algebraic
$\C^\times$-action. It follows from their results  that, when $\B$ is
coloop-free, one can always choose such a
circle action, which makes $Y(A,\theta)$, what we call, a {\em
  hyper-compact} hyperk\"ahler manifold. It means that the $U(1)$-action is Hamiltonian
with proper moment map with 
a minimum, and also that the holomorphic symplectic form $\omega_\C$
is of homogeneity $1$, meaning that for $\lambda\in \C^\times$
\be\label{homogeneity}\lambda^*\omega=\lambda \omega. \ee

For further results about the topology and geometry of toric
hyperk\"ahler varieties consult the papers \cite{bielawski-dancer},\cite{harada-proudfoot},
\cite{hausel-sturmfels}, \cite{hausel-swartz}  and \cite{konno}.

\section{Injective Hard Lefschetz for hyperk\"ahler manifolds}
\label{injective}

We are now ready to give two proofs of the following

\begin{theorem}For a smooth toric hyperk\"ahler variety $Y(A,\theta)$ of real dimension $4n-4d=4k$, 
such that $\B$ is coloop-free, we have that \be \begin{array}{c} L^{k-2i}:H^{2i}(Y,\C)\to H^{2k-2i}(Y,\C)\\ L^{k-2i}(\alpha)=\alpha\wedge\omega^{k-2i}\end{array}
\label{lefschetz}\ee is injective if $2i<k$, where $\omega=[\omega_I]$ is the cohomology class of the K\"ahler form corresponding to the complex structure $I$.  
\label{ghk}
\end{theorem}

Just like in Stanley's proof of the McMullen conjecture, we also have the following numerical consequences:

\begin{corollary}
 The $h$-vector $(h_1(\M),\dots,h_k(\M))$ of a coloop-free and rank
 $k$ matroid $\M$, which is (unimodularly and)
 rationally representable,
 satisfies \be\label{less}h_i(\M) \leq h_j(\M),\ee for $i\leq j \leq k-i$ and  
the $g$-inequalities (\ref{ginequalities}). 
\end{corollary}

\paragraph{\it Proof of Corollary.} Let the (unimodular) vector configuration
$\B=\{b_1,\dots,b_n\}\in
 \Z^k\subset \Q^k $ represent the matroid $\M$. Choosing a Gale dual
 configuration $\A=(a_1,\dots,a_n)\in \Z^d$ and a generic $\theta\in
 \N\A$, we can construct a smooth toric hyperk\"ahler variety $Y(A,\theta)$,
 whose Betti numbers agree with the $h$-numbers of $\M$. Now
 Theorem~\ref{ghk} immediately implies (\ref{less}). 
From Theorem~\ref{ghk} we can also deduce (\ref{ginequalities})
 exactly as in Stanley's argument for simplicial convex polytopes. See
 the introduction or for more details \cite[Theorem III.1.1]{stanley-book}.  
\endproof

\paragraph{\it Proof 1 of Theorem~\ref{ghk}.} As explained above we
 have a $\C^\times$-action on $Y:=Y(A,\theta)$, for which the
 corresponding $U(1)\subset \C^\times$-action is hyper-compact. Recall
 that this means that it is Hamiltonian with a proper moment $\mu_\R: Y\to \R$  
map with respect to $\omega$, and for which
the holomorphic symplectic form $\omega_\C$ is of homogeneity $1$ meaning (\ref{homogeneity}). Suppose that the fixed point set of the circle action has $f$ components, which are denoted by
$F_1,\dots, F_f$. The numbering is such that $\mu_\R(F_m)>\mu_\R(F_l)$ implies $m>l$. Now we define the Bialynicki-Birula stratification of $Y$ with respect to our $\C^\times$-action. Namely
define $U_m=\{p\in Y | \lim_{\lambda\rightarrow 0}\lambda p\in F_m\}$, which is an affine bundle over $F_m$. Moreover we let $U_{\leq m}=\cup_{j\leq m} U_j$ and $U_{<m}=\cup_{j<m} U_j$, which are 
open subvarieties of $Y$. Because the moment map $\mu_\R$ is proper it follows that
$U_{\leq f}=Y$, i.e. that we get this way a stratification of $Y$. Finally we denote by $N_m$ the negative normal bundle of $F_m$. Because the holomorphic symplectic
form is of homogeneity $1$ with respect to our $\C^\times$-action,
 it follows (cf. \cite[Proposition 7.1]{nakajima-hilbert}) that \be\mbox{rank}_\C (N_m)+\dim_\C(F_m)={\frac{1}{2}} \dim_\C Y=k
\label{half}.\ee By induction on $m$ we prove that the map $L^{k-2i}$ in (\ref{lefschetz}) when restricted to $U_{\leq m}$ is injective for $2i<k$. For $m=1$ the statement is clear because by (\ref{homogeneity})
$U_1=T^*F_1$ thus $\dim_\C(F_1)=k$ and the statement follows from the traditional Hard Lefschetz theorem for the compact K\"ahler manifold $F_1$. Now suppose we have the required injectivity
of the map $L^{k-2i}$ on $U_{<m}$. Then consider the decomposition
$U_{\leq m}=U_{<m}\cup U_m$. 
From this decomposition, using the Thom isomorphism \be\label{thomiso}H^{2i}(U_{\leq
  m}, U_{<m};\C)\cong H^{2i-2n_m}(U_m,\C),\ee
we get the cohomology exact sequence:
 $$0 \rightarrow H^{2i-2n_m}(U_m,\C)\stackrel{\tau}{\rightarrow} H^{2i}(U_{\leq m},\C)\stackrel{r}{\rightarrow} H^{2i}(U_{<m},\C)\rightarrow 0,$$ where $n_m=\mbox{rank}_\C (N_m)$, $\tau$ is the 
Gysin map and $r$ is the natural restriction map on cohomology. Now suppose $2i<k$ and  $0\neq \alpha\in H^{2i}(U_{\leq m},\C)$. If $r(\alpha)\neq 0$, then by induction we can deduce that 
$L^{k-2i}(\alpha)\neq 0$. However if $r(\alpha)=0$, then there is a $\beta\in H^{2i-2n_m}(U_m,\C)$ such that $\tau(\beta)=\alpha$. However $U_m$ is homotopy equivalent with 
the smooth compact K\"ahler manifold $F_m$ and $\omega|_{F_m}$ is a
K\"ahler class. If we denote $f_m=\dim_\C F_m$,
then the Hard Lefschetz theorem for $F_m$  yields that 
$0\neq \beta\wedge \omega^{f_m-2(i-n_m)}=\beta\wedge
\omega^{k-2i+n_m}|_{F_m},$ because $f_m+n_m=k$ by (\ref{half}). Because $\tau$ is injective we get that $\tau(\beta\wedge \omega^{k-2i}|_{F_m})=\alpha\wedge \omega^{k-2i}|_{U_{\leq m}}\neq 0$. 

The result follows. \endproof

\begin{corollary} For a hyper-compact hyperk\"ahler manifold $M$
(such as e.g. toric hyperk\"ahler varieties or Nakajima's quiver varieties \cite{nakajima} or moduli spaces of Higgs bundles \cite{hitchin}) 
we have that $$\begin{array}{c}
L^{k-2i} : H^{2i}(M,\C)\to H^{2k-2i}(M,\C)\\ L^{k-2i}(\alpha)=\alpha\wedge\omega^{k-2i}\end{array}$$ is injective if $2i<k$, where $\omega=[\omega_I]$ is the class of the  K\"ahler form corresponding to the complex structure $I$.  
\end{corollary}

We now recall our original proof of Theorem~\ref{ghk} from
\cite[Theorem 7.4]{hausel-sturmfels} in the smooth case because we will use the idea in the final section.

\paragraph{\it Proof 2 of Theorem~\ref{ghk}.}  Let $X_1, \dots X_r$
denote the irreducible components of the core of $Y$. Let $\phi_i:H^*(Y,\C)\to H^*(X_i,\C)$ denote the natural restrictions. 
The heart of the proof of \cite[Theorem 7.4]{hausel-sturmfels} is then  that
\begin{equation}
\label{crucialclaim}
{\rm ker}(\phi_1) \, \cap \,
{\rm ker}(\phi_2) \, \cap \,\, \ldots \,\,
\, \cap \, {\rm ker}(\phi_r) \quad = \quad \{0\}.
\end{equation}
In \cite{hausel-sturmfels} we presented two proofs of this fact. One \cite[Proposition 3.4]{hausel-sturmfels} was a more general result for
semi-projective toric orbifolds and the proof goes similarly to our first Proof 1 of Theorem~\ref{ghk} above, i.e. uses Morse theory type considerations with induction.
It turns out that  \cite[Proposition 3.4]{hausel-sturmfels} is
equivalent with the fact that the bounded complex of the polytope (or
in our case the bounded complex of the affine hyperplane arrangement
${\mathcal H}(B,\psi)$) is always contractible. 
The second proof was given after equation (34) of \cite{hausel-sturmfels}, which showed that (\ref{crucialclaim}) is in fact equivalent with 
Stanley's result  \cite[Proposition III.3.2]{stanley-book} that the Stanley-Raisner ring of a matroid is level. 

Now we proceed as follows. For $2i< k$ take $\alpha\in
H^{2i}(Y,\C)$. Then because of (\ref{crucialclaim}) we have a $j$ so
that  $\phi_j(\alpha)\in H^{2i}(X_j,\C)$ is nonzero. 
But the traditional hard Lefschetz  theorem for the smooth compact
K\"ahler manifold $X_j$ implies that $\phi_j(\alpha \wedge \omega^{k-2i})\neq 0$.

The result follows. \endproof

\begin{remark} 1. \cite[Theorem 7.4]{hausel-sturmfels} proves the same result, in the way sketched above, for a rationally representable matroid, i.e. for toric hyperk\"ahler orbifolds, not just
for smooth toric hyperk\"ahler varieties. Here we restricted our
attention to the smooth case, because the other Proof 1 only works in this case. The reason is that (\ref{half}) could be false
in the orbifold case.

2. Proof 1 works for any hyper-compact hyperk\"ahler manifold, however an extension of Proof 2 in the general case is not immediate. Indeed the equivalent of (\ref{crucialclaim}) (perhaps in intersection cohomology) is not known for a general hyper-compact hyperk\"ahler manifold. 

3. Another consequence of (\ref{crucialclaim}), explained in \cite[Section 7]{hausel-sturmfels}, is that one can present the cohomology ring of $Y$, in terms of cogenerator polynomials
corresponding to the $X_i$, the components of the core. Indeed this
algebraic presentation is rather similar to a presentation of a pure
$O$-sequence, the only difference will be that we replace the
cogenerator polynomials by monomials.  This similarity will lead to the proof of 
Theorem~\ref{gstanley} below.      
\end{remark}

\section{Proof of the $g$-inequalities for matroid complexes}

\label{g} In this section we will use the geometrical idea from our first proof of Theorem~\ref{ghk} to prove the following generalization:

\begin{theorem} The $h$-vector $(h_1(\M),\dots,h_k(\M))$ of a coloop-free and rank $k$ matroid $\M$ satisfies (\ref{less}) and the $g$-inequalities (\ref{ginequalities}). \label{gswartz}
\end{theorem}

\begin{remark} This was first proven by Swartz \cite{swartz}, by using an algebraic version of Chari's \cite{chari} decomposition theorem of matroids. Here we will show, that \cite{chari} is in fact 
gives us the geometrical structure for a general matroid so that we can repeat our Morse theory type first proof of Theorem~\ref{ghk}. In fact this proof is similar to
Swartz's original proof. 
\end{remark}

\begin{proof} So let us first recall Chari's result
  \cite[Theorem 3]{chari}:
\begin{theorem}[Chari] A coloop-free matroid complex is PS-decomposable.
\end{theorem}
A pure rank-$k$ simplicial complex $\Sigma$ on a vertex set
$\{1,\dots,n\}$ 
is {\em PS-decomposable} if it 
can be covered by pure rank-$k$ simplicial subcomplexes 
$\Sigma=\cup^m_{i=1}\Sigma_i$, such that 
\begin{itemize}
\item $\Sigma_1$ is the poset-theoretic
product of boundaries of simplices (a 
PS-$k$-sphere in the terminology of \cite{chari}), while for each
$i=2,\dots,m$, $\Sigma_i$ is the poset-theoretic product  of a simplex and a 
PS-sphere (called a PS-ball in \cite{chari}), and
\item For $i\geq 2$, $\Sigma_i\cap\left(\cup_{j=1}^{i-1}\Sigma_j\right) =\partial \Sigma_i$, where $\partial\Sigma_i$ denotes 
the pure rank-$(k-1)$ simplicial complex (which is just a PS-sphere in this case) whose facets are the rank-$(k-1)$ faces of $\Sigma_i$ that are contained in only one facet of $\Sigma_i$.  
\end{itemize}
We will show that Theorem~\ref{gswartz} holds for PS-decomposable simplicial complexes, a result which was also mentioned by Swartz in \cite{swartz}. We will see that this PS-ear-decomposition is in fact a very good combinatorial analogue of the Morse stratification of $Y$ (or rather its Lagrangian core) used in Proof 1 of Theorem~\ref{ghk}.  

 We first make a 
\begin{definition} Let $R$ be a ring and $M$ be a graded
  $R$-module. Then we say that $M$ satisfies injective hard Lefschetz
  (IHL for short) around degree $k/2$ for $\omega\in R_1$ if the map
  \begin{eqnarray*}\begin{array}{c} L^{k-2i}:M_{i}\to M_{k-i} \\
  L^{k-2i}(\alpha)=\alpha\omega^{k-2i}\end{array}\end{eqnarray*}
is injective for $0<i\leq k/2$. 
\end{definition}

We will proceed by induction on $m$ to show that 

\begin{eqnarray}
\begin{array}{c}\mbox{ there is an 
l.s.o.p $(\theta_1,\dots,\theta_k)$ so that the graded 
ring
$\C[\Sigma]/(\theta)$ }   \\ \mbox{ satisfies  IHL around
  $k/2$ with 
$\omega=\sum_i {x_i}$.}
\label{induction}
\end{array}
\end{eqnarray}
 When $m=1$, then $\Sigma$ is just a poset-theoretic product of 
boundaries of simplices and, therefore $\C[\Sigma]$ 
can be thought of 
as the torus equivariant cohomology ring of a product of projective
spaces, while an l.s.o.p. $(\theta)$ can be chosen so that $\C[\Sigma]/(\theta)$ is just the cohomology ring of the product of projective spaces, then
$\omega=\sum x_i$ is just a K\"ahler class, so the classical Hard Lefschetz theorem proves (\ref{induction}).

Now suppose we know our statement for $m-1$ and consider a pure 
rank-{$k$} simplicial complex with a PS-ear-decomposition. Let us 
denote $\Sigma_{<m}= \cup_{j=1}^{m-1} \Sigma_j$. Consider the natural
surjective map  $\C[\Sigma]\to\C[\Sigma_{<m}]$. We think of the kernel
of this map as a graded $\C[x_1,\dots,x_n]$-module and denote it by $\C[\Sigma,\Sigma_{<m}]$. So we have the following exact sequence of graded $\C[x_1,\dots,x_n]$-modules:
$$0\to \C[\Sigma,\Sigma_{<m}] \to  \C[\Sigma]
\to\C[\Sigma_{<m}] \to 0.$$

We now claim that we can find an l.s.o.p
$(\theta)=(\theta_1,\dots,\theta_k)$ for  $\C[\Sigma]$ such that  in
both graded $\C[x_1,\dots,x_n]$-modules
$\C[\Sigma_{<m}]/(\theta)$ and $\C[\Sigma, \Sigma_{<m}]/(\theta)$ the
IHL for $\omega$ is satisfied around degree $k/2$ .

By induction we know that the set of $(\theta)$ which is an
l.s.o.p. for $\C[\Sigma_{<m}]$ and  
$\C[\Sigma, \Sigma_{<m}]/(\theta)$ satisfies IHL for $\omega$ 
is non-empty 
and clearly Zariski open in $\C^{nk}$. 
Consequently the set of $(\theta)$ which is an l.s.o.p for
$\C[\Sigma]$ and  
$\C[\Sigma_{<m}]/(\theta)$ satisfies IHL for $\omega$ 
is non-empty and Zariski open. It is also
clear that the set of $(\theta)$ which is an l.s.o.p for $\C[\Sigma]$
and    $\C[\Sigma, \Sigma_{<m}]/(\theta)$ satisfies IHL around degree
$k/2$ for $\omega$
is Zariski open. We now prove that it is in fact non-empty. Take the
natural map $\C[\Sigma_m]\to\C[\partial\Sigma_m]$ and denote by 
$\C[\Sigma_m,\partial\Sigma_m]$ the kernel. We think of this
kernel as an $\C[x_1,\dots,x_n]$-module by letting the variables $x_j$
which correspond to vertices not in $\Sigma_m$ acting trivially. Then
it is easy to see that $\C[\Sigma_m,\partial\Sigma_m]$ and
$\C[\Sigma,\Sigma_{< m}]$ are isomorphic as graded
$\C[x_1,\dots,x_n]$-modules (this is the analogue of excision in
cohomology). But $\Sigma_m=\Delta\times \Phi$ 
is a poset-theoretic product of a $k$-simplex $\Delta$ with a
poset-theoretic product of boundary of 
simplices $\Phi$. Now it is clear that  
$$\C[\Sigma_m,\partial\Sigma_m]\cong \C[\Phi]\otimes\C[\Delta,\partial\Delta]$$ as graded 
$\C[x_1,\dots,x_n]$-modules (this corresponds to the Thom isomorphism (\ref{thomiso})
in cohomology). If $x_1,\dots, x_l$ correspond to the vertices of
$\Delta$ then 
$\C[\Delta,\partial\Delta]$ is just a free $\C[x_1,\dots,x_l]$-module
  generated by a degree $k$ element $x_1x_2\dots x_l$ (which is the
  analogue of the Thom class). 

First we note that the set of
  $(\theta)=(\theta_1,\dots,\theta_k)\in (\C[\Sigma])_1^k=\C^{nk}$ for
  which
  $$\C[\Sigma,\Sigma_{<m}]/(\theta):=\C[\Sigma,\Sigma_{<m}]/(\theta_1
  \C[\Sigma,\Sigma_{<m}]+\dots +\theta_k\C[\Sigma,\Sigma_{<m}])$$ 
satisfies IHL around degree $k/2$ for $\omega=\sum_{i=1}^n x_i$ is clearly Zariski
  open in $\C^{nk}$. Now we show that it is non-empty. Take 
$(\theta)=(x_1,\dots,x_l,\theta_{l+1},\dots,\theta_k)$, so that 
$ (\theta_{l+1},\dots,\theta_k)$ is an $l.s.o.p$ for $\C[\Phi]$ and
  $\C[\Phi]/(\theta_{l+1},\dots,\theta_k)$ satisfies IHL around
  $(k-l)/2$ with $\omega=\sum x_i $.  

 For this choice we have
  $$\C[\Sigma,\Sigma_{<m}]/(\theta)= x_1x_2\dots x_l
  \C[\Phi]/(\theta_{l+1},\dots,\theta_k),$$ and so IHL 
for $ \C[\Phi]/(\theta_{l+1},\dots,\theta_k)$ around degree $(k-l)/2$ 
implies IHL 
  for  $\C[\Sigma,\Sigma_{<m}]/(\theta)$ around degree $k/2$ with
  $\omega=\sum x_i$.

As the intersection of non-empty Zariski subsets of $\C^{nk}$ are non-empty 
we can choose a $(\theta)=(\theta_1,\dots,\theta_k)$, which is an
l.s.o.p for $\C[\Sigma]$ and $\C[\Sigma_{<m}]$ and for which both
$\C[\Sigma_{<m}]/(\theta)$ and $\C[\Sigma,\Sigma_{<m}]/(\theta)$
satisfies IHL around $k/2$ with $\omega=\sum x_i$. Now using the short exact
sequence: 

$$0\to \C[\Sigma,\Sigma_{<m}]/(\theta) \to  \C[\Sigma]/(\theta)
\to\C[\Sigma_{<m}]/(\theta) \to 0,$$ we can repeat the argument of
Proof 1 of Theorem~\ref{ghk}, to get that $\C[\Sigma]/(\theta)$
satisfies IHL around $k/2$ with $\omega=\sum x_i$. 

Because a PS-decomposable simplicial complex $\Sigma$ is shellable (see
\cite[Proposition 5]{chari}), and
so Cohen-Macaulay, we have that
$h_i(\Sigma)=\dim_\C( (\C[\Sigma]/(\theta))_i)$ and so Theorem~\ref{gswartz}
follows. \end{proof}

\begin{remark} Because we have Hard Lefschetz theorem for boundary complexes of
  simplicial convex polytopes the above proof would have worked
  equally well for simplicial complexes with a decomposition just like
PS-ear-decomposition above, but changing $PS$-spheres, in the definition, with
  boundary complexes of simplical convex polytopes. For a unimodularly
  and rationally
  representable matroid such a presentation always arises naturally. Namely we
  can 
  consider the Morse stratification (for details on this see
  \cite{harada-proudfoot}) of a hyper-compact $U(1)$-action on 
the bounded complex of a generic hyperplane arrangement,
  representing our given matroid. In this case the above combinatorial
  proof of Theorem~\ref{gswartz} would essentially agree with Proof 1
  of Theorem~\ref{ghk}. 
\end{remark}

\section{Proof of the $g$-inequalities for pure $O$-sequences} 
\label{stanley}

First a definition:

\begin{definition} A sequence of non-negative integers $(h_1,h_2,\dots, h_k)$ is called a pure $O$-sequence, if $h_k>0$ and there exists monomials $m_1,\dots,m_{h_k}$ of
degree $k$ in the degree one variables $x_1,\dots,x_{h_1}$, so that 
\begin{multline*} h_l=\#\{ m | m \ \mbox{is a monomial of degree $l$ in variables } x_1,\dots,x_{h_1},\\ \mbox{ such that } 
m|m_i, \mbox{ for some } 0<i\leq h_k.\}
\end{multline*}
\end{definition}

Now we can state a long standing conjecture of Stanley \cite{stanley-conjecture}:

\begin{conjecture}[Stanley] The $h$-vector $(h_1(\M),\dots,h_k(\M))$ of a rank $k$ matroid $\M$ is a pure $O$-sequence. 
\end{conjecture}

This conjecture is still open for general matroids, though recently it
has been deduced for cographic matroids using \cite{lopez}, i.e. for the Betti numbers
of toric quiver varieties \cite[Section 8]{hausel-sturmfels}. Another attack on Stanley's conjecture has been to deduce numerical 
inequalities between the numbers in a pure $O$-sequence and then prove these inequalities for the $h$-vector of a matroid complex. 
 As an example Hibi \cite{hibi} proved that for a pure $O$-sequence one has \be\label{lessa} h_i\leq h_{j}, \ee where $i\leq j\leq k-i$ and in particular that  
$$h_1\leq h_2\leq \dots \leq h_{\lfloor{\frac{k}{2}}\rfloor},$$ this was in turn proven for $h$-vectors of matroid complexes by Chari \cite{chari}.

Here we strengthen this result by proving the following

\begin{theorem}  A pure $O$-sequence $(h_1,h_2,\dots, h_k)$ satisfies (\ref{lessa}) and the $g$-inequalities \label{gstanley}. 
\end{theorem}

\begin{corollary} 
Theorem~\ref{gswartz} is a consequence of Stanley's Conjecture~\ref{stanley}.
\end{corollary}

\paragraph{\it{Proof of Theorem~\ref{gstanley}.}} We are going to follow the
  structure of Proof 2 of Theorem~\ref{ghk}. Namely take a pure
  $O$-sequence  $(h_1,h_2,\dots, h_k)$ with generating monomials
  $m_1,\dots,m_{h_k}$ in variables $x_1,\dots, x_{h_1}$. First we
  construct a graded ring \bes \begin{array}{c} R={\frac{\C[\partial_1,\dots,
  \partial_{h_1}]}{I}  }\\ I= {\rm ann}(m_1)\cap\dots \cap {\rm ann}(m_{h_k})\end{array}\ees which will be the analogue of the cohomology ring $H^*(Y,\C)$ of a toric hyperk\"ahler manifold. Here $\partial_i$ is a variable of degree one, which we think of as a differential operator, satisfying
$\partial_i(x_j)=\delta_{ij}$. The ideal in the denominator is the ideal $I$ of polynomials in the $\partial_i$ which annihilate all the monomials $m_j$. Clearly $\dim R_j=h_j.$ 
Then we construct graded rings \begin{eqnarray*}\begin{array}{c}R^j={\frac{\C[\partial_1,\dots, \partial_{h_1}]}{I_j} }\\ I_j= \mbox{ann}(m_j)\end{array}\end{eqnarray*} for each monomial $m_j$, which will be the analogue of $H^*(X_j,\C)$ (in fact it is useful to think about $R^j$ as the cohomology ring of the product of projective spaces of dimension given by the 
exponents in the monomial $m_j$). Because $I\subset I_j$, we have a natural map $p_j:R\to R^j$. The equation $I=\cap_j I_j$ now implies the analogue of (\ref{injective}), i.e. that  the map $$p=p_1\times\dots \times p_{h_k}: R\to R^1\times\dots \times R^{h_k}$$ is injective. 
Now take the degree $1$ class $\omega=\sum_j \partial_j$. It is clear
  that the map 
$L^{k-2i}_j:R^j_i\to R^j_{k-i} $ given by  $L^{k-2i}_j(\alpha)=\alpha  p_j(\omega^{k-2i})$ is injective for $2i<k$. 
It follows e.g. if we think
of $R^j$ as the cohomology ring of the product of projective spaces. 
Then $p_j(\omega)$ corresponds to
the natural ample class, so the hard Lefschetz theorem implies
  injectivity of $L^{k-2i}_j$. Of course in this case one can check this
  result by hand for the explicitly defined rings $R^j$.  The
  injectivity of $p$ and of $L^{k-2i}_j$ implies the injectivity of $L^{k-2i}:R_i
  \to R_{k-i}$, $L^{k-2i}(\alpha)=\alpha \omega^{k-2i}  $   
for $2i<k$. The result follows. \endproof


\begin{thebibliography}{1}

\bibitem{bielawski-dancer}
R.~Bielawski, A.~Dancer, {The geometry and topology of toric
 hyperk\"ahler manifolds},{\em Comm.  Anal. Geom.} {\bf 8} (2000), no. 4, 727-760.

\bibitem{chari} M.~Chari,
{Two decompositions in topological combinatorics with applications to matroid complexes.} 
\emph{Trans. Amer. Math. Soc}. {\bf 349} (1997), no. 10, 3925--3943.



\bibitem{griffiths-harris} P. Griffiths, J. Harris. {\em Principles
of algebraic geometry}. New York, Wiley 1978

\bibitem{harada-proudfoot} M.~Harada, N.~Proudfoot: Properties of the residual circle action on a toric hyperkahler variety,
(to appear in {\em Pac. J. Math.}), arXiv: math.DG/0207012 


\bibitem{fulton} W.~Fulton: {\em Introduction to Toric Varieties},
Princeton University Press, 1993.



\bibitem{hausel-sturmfels} T.~Hausel, B.~Sturmfels: Toric hyperk\"ahler varieties, {\em Documenta Mathematica}, {\bf 7} (2002), 495--534, arXiv:  math.AG/0203096 

\bibitem{hausel-swartz} T.~Hausel, E.~Swartz: Intersection forms of toric hyperk\"ahler varieties, 
preprint,\newline arXiv:math.AG/0306369 

\bibitem{hibi} T.~Hibi:  What can be said about pure $O$-sequences?  
  {\em J. Combin. Theory Ser. A}  {\bf 50}  (1989),  no. 2, 319--322.

\bibitem{hitchin} N. Hitchin. The self-duality equations on a Riemann 
surface, {\em Proc. London Math. Soc.} (3) {\bf 55} (1987) 59-126

\bibitem{konno} H.~Konno: Cohomology rings of toric hyperk\"ahler manifolds,
{\em Internat. J. Math.} {\bf 11} (2000), no. 8, 1001--1026.

\bibitem{lopez}
C.M.~Lopez: Chip firing and the Tutte polynomial,
{\em Ann.~Combinatorics} {\bf 1} (1997) 253--259.

\bibitem{mcmullen}
P.~McMullen, The numbers of faces of simplicial polytopes, {\em Israel J.
  Math.} \textbf{9} (1971), 559--570.

\bibitem{nakajima-hilbert} H.~Nakajima. {\em Lectures on Hilbert schemes of points
on surfaces.} University Lecture Series, 18. American Mathematical Society, Providence, RI,  1999.

\bibitem{nakajima} H.~Nakajima,
Quiver varieties and finite-dimensional representations of quantum affine algebras.\emph{J. Amer. Math. Soc.} {\bf 14} (2001), no. 1, 145--238



\bibitem{stanley-toric}
R.~Stanley, The number of faces of a simplicial convex polytope, {\em Adv. in
  Math.} \textbf{35} (1980), no.~3, 236--238.

\bibitem{stanley-conjecture}
R.~Stanley, Cohen-{M}acaulay complexes, in \emph{Higher combinatorics} (Proc.
  NATO Advanced Study Inst., Berlin, 1976), Reidel, Dordrecht, 1977,
  pp.~51--62. NATO Adv. Study Inst. Ser., Ser. C: Math. and Phys. Sci., 31.

\bibitem{stanley-book} R.P. Stanley:
{\em Combinatorics and Commutative Algebra},
2nd ed., Birkh\"auser Boston, 1996.

\bibitem{swartz} E.~Swartz, $g$-elements of matroid complexes, (to appear in {\em J. Comb. Theory Ser. B})


\end{thebibliography}
\end{document}